\documentclass[journal,twoside,web]{ieeecolor}
\usepackage{generic}
\usepackage{cite}
\usepackage{amsmath,amssymb,amsfonts}
\usepackage{algcompatible}
\usepackage{algorithm}
\usepackage{graphicx}
\usepackage{hyperref}
\hypersetup{hidelinks=true}
\usepackage{textcomp}
\usepackage{bm}
\usepackage{adjustbox}
\usepackage{amssymb}
\definecolor{Red}{rgb}{1,0,0}
\definecolor{Green}{rgb}{0,.6,0}
\definecolor{Blue}{rgb}{0,0,0}

\newcommand{\norm}[1]{\left\lVert#1\right\rVert}
\newcommand{\normy}[1]{\left\lVert#1\right\rVert_{\mathbf{D}(\yy)}}

\newtheorem{theorem}{\bf{Theorem}}
\newtheorem{lemma}{\bf{Lemma}}

\DeclareMathOperator*{\argmin}{arg\,min}

\newcommand{\RR}{\mathbb{R}}

\newcommand{\xx}{\mathbf{x}}
\newcommand{\bb}{\mathbf{b}}
\newcommand{\cc}{\mathbf{c}}

\newcommand{\zz}{\mathbf{z}}
\renewcommand{\AA}{\mathbf{A}}
\newcommand{\HH}{\mathbf{H}}
\newcommand{\Df}{\textcolor{Blue}{\mathcal{D}_t}}

\newcommand{\vv}{\bm{\nu}}

\newcommand{\yy}{\mathbf{y}}

\newcommand{\DD}{\mathbf{D}}

\newcommand{\nn}{\eta}
\newcommand{\WW}{\mathbf{W}}

\def\BibTeX{{\rm B\kern-.05em{\sc i\kern-.025em b}\kern-.08em
    T\kern-.1667em\lower.7ex\hbox{E}\kern-.125emX}}
\begin{document}

\title{Online Interior-point Methods for Time-Varying Equality-constrained Optimization}

\author{Jean-Luc Lupien, \IEEEmembership{Student Member, IEEE}, Iman Shames \IEEEmembership{Member, IEEE}, and Antoine Lesage-Landry, \IEEEmembership{Member, IEEE}
\thanks{Submitted July 2023. This work was funded by NSERC and IVADO.}
\thanks{J-L. Lupien and A. Lesage-Landry are with the Department of Electrical Engineering, Polytechnique Montréal, MILA \& GERAD, Montréal, QC, Canada, H3T 1J4. e-mail: \texttt{jean-luc.lupien@polymtl.ca}, \texttt{antoine.lesage-landry@polymtl.ca}}
\thanks{I. Shames is with CIICADA Lab in the School of Engineering at the Australian National University, Canberra, ACT, Australia. e-mail: \texttt{iman.shames@anu.edu.au}}}

\maketitle

\begin{abstract}
An important challenge in the online convex optimization (OCO) setting is to incorporate generalized inequalities and time-varying constraints. The inclusion of constraints in OCO widens the applicability of such algorithms to dynamic and safety-critical settings such as the online optimal power flow (OPF) problem. In this work, we propose the first projection-free OCO algorithm admitting time-varying linear constraints and convex generalized inequalities: the online interior-point method for time-varying equality constraints (\texttt{OIPM-TEC}). We derive simultaneous sublinear dynamic regret and constraint violation bounds for \texttt{OIPM-TEC} under standard assumptions. For applications where a given tolerance around optima is accepted, we \textcolor{Blue}{employ an alternative} OCO performance metric -- the epsilon-regret -- and a more computationally efficient algorithm, the epsilon-\texttt{OIPM-TEC}, that possesses sublinear bounds under this metric. Finally, we showcase the performance of these two algorithms on an online OPF problem and compare them to another OCO algorithm from the literature.
\end{abstract}


\section{Introduction}
\label{sec:introduction}

\IEEEPARstart{O}{nline} convex optimization (OCO) with time-varying constraints, henceforth referred to as constrained OCO, is an emerging field with many practical applications including network resource allocation~\cite{NetworkExample}, portfolio selection~\cite{hazanOCO}, and control of electric grids~\cite{OnlinePowerReview} or of multi-energy systems~\cite{L1penalty}. In OCO, an agent makes sequential decisions to minimize an environment-determined loss function~\cite{hazanOCO, Shalev, zinkevich2003online}. Constrained OCO adds the requirement that the decision sequence is subject to time-varying constraints. Integrating online or time-varying constraints is essential to ensure effective implementation of OCO algorithms in safety-critical applications. The performance of a constrained OCO algorithm is measured both in terms of its regret and its constraint violation~\cite{DynamicPolyak, Shalev}. These represent the cumulative distance from optimality and from feasibility, respectively. Regret can be calculated either by comparing the decision sequence to the best fixed decision in hindsight (static regret) or to the round-optimal decision (dynamic regret)~\cite{hall2015online}. The dynamic regret is a tighter metric than the static regret as a bounded dynamic regret guarantees a bounded static regret~\cite{MOSP}. When designing an OCO algorithm, the objective is to obtain sublinear regret and constraint violation bounds which imply the time-averaged regret and constraint violation will go to zero as the time-horizon goes to infinity~\cite{hazanOCO, Shalev, MOSP}. This means implemented decisions will be optimal and feasible, on average, over a long time horizon.

A promising application of constrained OCO is the online optimal power flow (OPF) problem. The OPF problem consists of finding the optimal power generation profile that minimizes costs for the operator of the electric grid while ensuring that all loads are adequately supplied and that physical and safety limits are respected~\cite{ConvexRelaxations}. This problem is non-convex and thus NP-hard to solve. For this reason, multiple convex relaxations of the OPF have been formulated to compute an approximated solution in polynomial time~\cite{ConvexRelaxations, InteriorforOPF, taylor2015convex}. However, the tightest convex relaxations are conic relaxations that leverage generalized inequalities~\cite{taylor2015convex}. This makes solving these problems tedious and time-consuming in large-scale systems. In modern power grids, the integration of renewable energy sources introduces power fluctuations on a significantly faster timescale when compared to conventional generators~\cite{taylor2016power}. This makes the traditional approach of solving the OPF every 15-30 minutes~\cite{OnlinePowerReview, ConicSurvey} unable to conform to the renewable timescale. The efficiency, adaptability, and resiliency of constrained OCO algorithms can overcome this limitation.

However, solving a conic relaxation of the OPF problem in an OCO context is challenging for multiple reasons. First and foremost is the integration of time-varying equality constraints linked to changes in power demand. Addressing the generalized inequalities associated with the conic constraints is also an obstacle. Many OCO algorithms use a projection step~\cite{MOSP, LagrangianTime-Varying, QueueConstrained} onto the feasible space defined by these constraints. This is problematic because when the feasible space is the cone of positive semi-definite (PSD) matrices or a second-order cone -- as is the case in the relaxed OPF -- the projection can be as costly as solving the offline problem. Additionally, the less stringent static regret is ill-adapted to an online OPF setting. A fixed solution satisfying the constraints at all timesteps may incur arbitrarily large losses or might not exist at all~\cite{MOSP}. Therefore, a candidate OCO algorithm for the online OPF is \emph{one admitting time-varying equality constraints, being projection-free, and possessing a sublinear dynamic regret and constraint violation}. We propose an online optimization approach based on a primal interior-point method that possesses these three qualities.

\paragraph*{Related work}

Although many constrained OCO approaches can be found in the literature, none possesses all three required properties necessary to effectively solve the online OPF. The majority of OCO algorithms that admit time-varying constraints use static regret as a performance metric. This metric uses the best fixed decision in hindsight as the comparator sequence. A common approach in this setting is to use a virtual queue to ensure a bounded static regret while satisfying long-term or slowly varying constraints. This approach is equivalent to solving a modified Lagrangian system. Virtual queues are used in~\cite{JMLR, FTRLTime-Varying, OnlineConvexBandit,OnlineEpsilonProof} to achieve sublinear static regret and constraint violation bounds. The model-based augmented Lagrangian method (\texttt{MALM}) presented in~\cite{LagrangianTime-Varying}, which possesses sublinear regret and constraint violation bounds, has been shown to experimentally outperform many of these previous methods. A projection-free OCO algorithm utilizing barrier functions as regularizers is presented in~\cite{HazanInterior} and achieves sublinear static regret. This approach does not admit time-varying equality constraints or a dynamic regret bound. 

Algorithms with dynamic regret bounds have also been developed to handle time-varying constraints. An online saddle-point method (\texttt{MOSP}) proposed in~\cite{MOSP} and an approach using virtual queues presented in~\cite{QueueConstrained} possess sublinear dynamic regret and constraint violation bounds. These bounds are, however, predicated on strict limits on objective and constraint function variations as well as time-varying stepsizes. Both these algorithms rely on a costly projection step for conic constraints. Sublinear dynamic regret and constraint violation is obtained in~\cite{DynamicPolyak} even for non-convex functions. The bounds are looser than those of~\cite{MOSP, QueueConstrained} and the algorithm requires a projection at every iteration. A primal-dual mirror descent approach was proposed in~\cite{DistributedOnline2} for distributed OCO and analyzed using dynamic regret but the authors only obtained constraint violation bounds that grow linearly with time. A second-order method was presented in~\cite{OnlineNewton} for unconstrained optimization which was extended in~\cite{OPENM} to admit time-varying equality constraints. Both these approaches have tight dynamic regret bounds but are incompatible with inequality constraints.

Approaches to iteratively solve the OPF have also been presented. An algorithm based on a quasi-Newton update was used to solve the OPF in real-time in~\cite{realtimeOPF}. Although this algorithm has good tracking performance under strict conditions, the authors did not derive a bound on constraint violation. A distributed feedback controller is designed in~\cite{OptimalPowerFlowPursuit} but uses a linear approximation of the OPF rather than a conic relaxation. A model-free online feedback optimization algorithm is presented in~\cite{OFO}. This approach is shown to converge to optimal solutions of the OPF problem but does not have theoretical guarantees on tracking or constraint violation.

\textcolor{Blue}{Finally, interior-point methods applied to systems with time-varying constraint and objective functions were presented in \cite{InteriorOnline} and~\cite{fazylab2}. However, this framework differs from OCO as continuous systems are considered where current-round information is available. Another method considering time-varying linear equality and inequality constraints  is proposed in~\cite{casti2024control}. It, however, does not provide a convergence proof for inequality-constrained problems and further uses a projection step.}

In this paper, we present the first constrained OCO algorithm capable of satisfying all three necessary conditions for effective resolution of the online OPF. This is achieved by extending the \texttt{OPEN-M} algorithm from~\cite{OPENM} to integrate generalized inequalities using an interior-point method-like update. Our contributions can be summarized as follows:

\begin{itemize}
	\item We present the first projection-free OCO algorithm that admits convex generalized inequalities and time-varying linear equality constraints while achieving tight sublinear dynamic regret and constraint violation bounds.
 \item We \textcolor{Blue}{employ} the notion of epsilon-regret in which one accepts being in an epsilon-neighbourhood of the optimal solution and propose an efficient algorithm possessing tight epsilon-regret and constraint violation bounds.
	\item We illustrate the performance of our approaches by solving an online OPF based on the 33-bus Matpower case.
	
\end{itemize}

\section{Preliminaries}

In this section, we define useful terms, lay down important assumptions, and introduce the Newton's step which will be used in the analysis of  our main results.
\subsection{Definitions}

\textcolor{Blue}{
In this work, our goal is to solve time-varying equality constrained OCO problems subject to generalized inequality constraints described in~\eqref{prob1}. The generalized inequality constraints can, for example, be defined with respect to the cone of positive semidefinite matrices or a second-order cone. Formally, let $\cc \in \RR^N$ be our objective vector with $N\in\mathbb{N}$ the dimension of the decision variable. Let $\xx_t \in \RR^N$  be the decision vector at any time $t = 0,1,2,...,T$ where $T\in\mathbb{N}$ is the time horizon. Let $K_i$, $i=1,2,...,I$, $I\in\mathbb{N}$ be a proper cone as defined in~\cite[Section 3.6.2]{boyd}. All functions $g_i(\xx) : \RR^N\mapsto \RR^{K_i}$, $i=1,2,...,I$, are assumed to be convex with respect to $K_i$ and self-concordant. Self-concordancy is defined in \eqref{def:selfconcordance} and is a common assumption in the interior-point literature~\cite{boyd, renegar}. Let the matrix $\AA \in \RR^{P\times N}$ be a full row-rank matrix with $P\in \mathbb{N}$, $P<N$. Let the vector $\bb_t \in \RR^P$ be time-varying. In the sequel, the norm $\norm{\cdot}$ denotes the Euclidean norm. In summary, we tackle the following online optimization problem:
}
\begin{align}\label{prob1}
\begin{split}
    \min_{\xx_t}&\quad \cc^\top\xx_t\\
    \text{s.t.} \quad & g_i(\xx_t)\preceq_{K_i} 0, \quad \text{for } i=1,2,3,...,I\\
    & \AA \xx_t - \bb_t = \mathbf{0}.
    \end{split}
\end{align}
\textcolor{Blue}{We denote the feasible set at time $t$ by $\Df \stackrel{\text{def}}{=} \left\{\xx \in \RR^N | g_i(\xx) \preceq 0, \forall i, \AA\xx-\bb_t=\mathbf{0}\right\}$. We remark that the only time-dependent parameter is the vector $\bb_t$ in this formulation. Although this limits the scope of OCO problems that can be represented, many applications are adequately modelled in this manner. This is the case for the optimal power flow problem, for example, where the only time-dependent parameter are the loads that need to be supplied~\cite{taylor2015convex, ConvexRelaxations}.}
\textcolor{Blue}{We define the set of optimal decisions, $\mathcal{X}_t\in\mathbb{R}^N$, to \eqref{prob1} at every timestep as:
\begin{align*}
\mathcal{X}_t \stackrel{\text{def}}{=} \left\{\xx_t^*\in\mathbb{R}^N| \xx_t^* \in \argmin \cc^\top \xx_t^*, \xx_t^* \in \Df\right\}.
\end{align*}}
%\textcolor{Blue}{This optimal decision might not be unique. (Problem?)}

%To measure the performance of constrained OCO algorithms, two metrics are employed: regret and constraint violation. These represent the cumulative distance from optimality and from feasibility, respectively. 
In this work, dynamic regret is used which takes the round-optimal solution $\xx_t^*$ as a comparator. This can be formally defined as:
\begin{equation}\label{def:regret}
    R_{\text{d}}(T) = \sum_{t=1}^T f(\xx_t)-f(\xx_t^*).
\end{equation}
\textcolor{Blue}{Let $d_{K_i}(\xx) = \min\limits_{\tilde{\xx}} \norm{\xx-\tilde{\xx}} \text{ s.t. } g_i(\tilde{\xx}) \preceq_{K_i} 0$, the closest distance from the point $\xx$ to a point in the interior of the proper cone $K_i$.
We can then define constraint violation as:
\begin{equation}
    \text{Vio}(T) = \sum_{t=1}^T\left[ \sum_id_{K_i}(\xx_t)+\norm{\AA\xx_t-\bb_t}\right].
\end{equation}
Constraint violation at $t$ is null if decision $\xx_t$ is feasible and positive otherwise.} \textcolor{Blue}{Useful for our analysis is the cumulative variation between optimal decisions $V_T$  which we define as:
\[V_T = \sum_{t=1}^{T} H_d(\mathcal{X}_t, \mathcal{X}_{t-1}),\]
where $H_d$ represents the Hausdorff distance between sets. Of note, when the optimal set is a singleton, we retrieve the more common definition of cumulative variation given by:
\[V_T = \sum_{t=1}^T \norm{\xx_t^*-\xx_{t+1}^*}.\]
When considering dynamic regret, a bound of the order of $V_T$ means that the algorithms successfully tracks the optimal solution through time \cite{hazanOCO, MOSP, OPENM}. It is often assumed $V_T$ is sublinear.
}
We define the cumulative variation in the right-hand term of the time-varying  equality constraint $\bb$ as: 
\[V_\bb = \sum_{t=1}^{T} \norm{\bb_t-\bb_{t-1}}.\]
\textcolor{Blue}{
Let $\norm{\cdot}_{\HH}$ denote the matrix norm where $\HH\in \RR^{N\times N}$ is a positive definite matrix. For any vector $\zz \in \RR^N$, the norm is defined as: \[\norm{\zz}_\HH = \sqrt{\zz^\top \HH \zz}.\] Notably, as in~\cite[Section 2.1]{renegar}, this norm is well-defined on convex domains for any twice differentiable strongly convex function $f(\xx):\RR^N\mapsto \RR$ by using the Hessian $\nabla^2 f(\xx)$. This so-called Hessian norm is: $\norm{\zz}_{\nabla^2 f(\xx)}^2 = \zz^\top \nabla^2 f(\xx)\zz$, and will be used extensively in the analysis of our proposed constrained OCO algorithms. }

An important property of functions that guarantees convergence of interior-point methods is self-concordance~\cite{boyd, renegar}. Self-concordance can be defined as all \textcolor{Blue}{differentiable} functionals $f$ for which any $\xx_1$, $\xx_2$ $\in\RR^N$ such that $\norm{\xx_2-\xx_1}_{\nabla^2 f(\xx_1)} < 1$, satisfy the following condition:
\begin{align}
\label{def:selfconcordance}
    \frac{\norm{\mathbf{v}}_{\nabla ^2 f(\xx_2)}}{\norm{\mathbf{v}}_{\nabla^2 f(\xx_1)}} \le \frac{1}{1-\norm{\xx_2-\xx_1}_{\nabla^2 f(\xx_1)}},
\end{align}
for any non-zero vector $\mathbf{v}\in\RR^N$. \textcolor{Blue}{This definition of self-concordance generalizes the more common definition provided in~\cite[Section 9.6]{boyd}. This is expressed as a function $f(x) : \RR\mapsto\RR$ being self-concordant if $|f'''(x)|\le 2f''(x)^{\frac{3}{2}}.$ A function $f(\xx)\in\RR^N\mapsto\RR$ is said to be self-concordant if it is self-concordant along any line in its domain.}
\textcolor{Blue}{Next, we define the barrier functional $\phi(\xx)$ such that $\text{dom } \phi(\xx)= \{\xx\in\RR^n | g_i(\xx) < 0\}$ as follows:}
\textcolor{Blue}{$\phi(\xx) = \sum_{i=1}^I -\log\big(-g_i(\xx)\big).$}
\textcolor{Blue}{Here the logarithm is either the natural logarithm if $g_i$ is a scalar function or a generalized logarithm as defined in~\cite[Section 11.6]{boyd}. This barrier functional is convex and self-concordant over its domain~\cite{renegar}}.

Interior-point methods do not attempt to solve \eqref{prob1} directly. Instead, a sequence of convex optimization problems related to the original problem are solved. This is achieved in the following way. Let $\nn\in \mathbb{R}^+$ be a positive barrier parameter. This barrier parameter is useful in defining the following equality-constrained problem:
\begin{align}\label{prob2}\begin{split}\min_{\xx_t}\quad &\nn\cc^\top\xx_t+\phi(\xx_t)\\
\text{s.t.} \quad& \AA \xx_t - \bb_t = \mathbf{0}.
\end{split}
\end{align} Problem \eqref{prob2} is of interest because when $\nn \rightarrow +\infty$, the solutions of \eqref{prob1} and \eqref{prob2} coincide. The strategy is, therefore, to incrementally increase $\nn$, iteratively solving the original problem.

\textcolor{Blue}{
In the sequel, we use the following notation:
\begin{equation}\label{defd}d(\xx_t, \nn) = d_\nn(\xx_t) \stackrel{\text{def}}{=} \nn \cc^\top\xx_t+\phi(\xx_t).\end{equation}
This represents the objective function we seek to minimize at every timestep $t$ given a barrier parameter $\nn$.}

\subsubsection{Assumptions}
Certain sufficient conditions need to be met to ensure the regret and constraint violation suffered by the algorithms are sublinear. These are standard in the interior-point literature~\cite{InteriorforOPF} and will hold for any convex relaxation of the OPF problem~\cite{ConvexRelaxations}.

\begin{enumerate}

    \item The relative interior of the feasible set $\Df$ forms a closed and convex set for all time.
    \item The feasible set $\Df$ is non-empty for all time. This implies there exists at least one feasible point and that Slater's condition holds.

    \item The barrier functional $\phi(\xx)$ is a self-concordant barrier functional of complexity $v_f\in\RR^+$. That is,
    \begin{align}\label{def:sccomplexity}\begin{split}
        \sup_\xx\quad \nabla\phi(\xx)\nabla^2\phi^{-1}(\xx)\nabla\phi(\xx)&\le v_f.
        \end{split}
    \end{align}

    The complexity of standard barrier functionals can be reasonably bounded by the dimension of the decision space, i.e., $v_f=O(N)$~\cite{boyd}. For example, the barrier functionals used for the space of positive semi-definite (PSD) matrices, second-order cones, quadratic functions, and linear functions all satisfy \eqref{def:sccomplexity}~\cite[Section 2.3]{renegar}. \textcolor{Blue}{Of note, this implies that on the set $\Df$, $\phi(\xx)$ is strongly convex~\cite{boyd, renegar}.}

    \item \textcolor{Blue}{For $\xx \in \Df$, we have: \[\norm{\begin{bmatrix} \nabla^2 \phi(\xx) & \AA^\top\\ \AA& 0\end{bmatrix}} \le \frac{1}{m},\]
    where $m>0$. This assumption is similar to~\cite[Chapter 10]{boyd} and~\cite{InteriorOnline}.}

\end{enumerate}

We assume that Assumptions 1-4 hold hereinafter.

\subsubsection{Newton Step}

Define vector $\yy\in \mathbb{R}^{N+P}$ as: $\yy = \begin{bmatrix}\xx\\ \vv\end{bmatrix}$. Let $\Gamma_\nn^t(\yy) \in \mathbb{R}^{N+P}\mapsto \mathbb{R}$, the Lagrangian function, be: \[\Gamma_\nn^t(\yy) = d_\nn(\xx)+\vv^\top (\AA\xx-\bb_t).\] This function has a saddle-point for $\xx=\xx_\nn^*,$ where $\xx_\nn^*$ is an optimum of \eqref{prob2}.

The Lagrangian function has the following interesting properties:
\begin{align*}
\nabla \Gamma_\nn^t(\yy)\stackrel{\text{def}}{=}r_t(\yy,\nn) &= \begin{bmatrix}
    \nabla d_\nn(\xx) + \AA^\top \vv\\
    \AA \xx -\bb_t
\end{bmatrix}\\
\nabla^2\Gamma_\nn^t(\yy)\stackrel{\text{def}}{=}\DD(\yy) &= \begin{bmatrix} \nabla^2 \phi(\xx) & \AA^\top\\ \AA& 0\end{bmatrix}.\end{align*}
We remark that $\DD(\yy)$ is the gradient of $r_t$ and is independent of both time $t$ and barrier parameter $\nn$. %The Hessian norms in the sequel will be defined with respect to this Lagrangian function, i.e. $\norm{\yy_1}^2_{\yy_2} = \yy_1^\top\DD(\yy_2)\yy_1,\quad\forall \yy_1,\yy_2\in \Df.$

Furthermore, the third derivative of $\Gamma_\nn^t(\yy)$ coincides with the
third derivative of $d_\nn (\xx)$ having no dependency in $\vv$. \textcolor{Blue}{Because
the third derivative of $d_\nn (\xx)$ is bounded on $\Df$ by Assumption 3, the third derivative of $\Gamma_\nn^t(\yy)$ is also bounded.} Therefore, we have by~\cite[Section 2.5]{renegar}, that $\Gamma_\nn^t(\yy)$ is self-concordant

The infeasible Newton's step, $n_t(\yy_t,\nn)$, -- so-called because it does not require a feasible initial point -- at round $t$ associated to a barrier parameter $\nn$ can now be defined as:
\[n_t(\yy_t,\nn)=\Delta \yy_t = \textcolor{Blue}{-\DD(\yy_t)^{-1}}r_t(\yy_t,\nn),\]
where $t$ denotes the current round.

We now present a sequence of lemmas that will later be used to to prove our main results in Section \ref{algos}.

%\begin{lemma}
%    The Lagrangian function
%\end{lemma}
%\begin{proof}
%    The third derivative of $\Gamma(\yy)$ coincides with the third derivative of $d_\nn(\xx)$ having no dependency in $\vv$. Because the third derivative of $d_\nn(\xx)$ is bounded by assumption, this provides a bound on the third derivative of $\Gamma(\yy)$. By~\cite[Section~2.5]{renegar}, we have that $\Gamma(\yy)$ is self-concordant.
%\end{proof}

\begin{lemma}\label{invHess}
    Let $\Gamma_\nn^t(\yy) = d_\nn(\xx)+\vv^\top(\AA\xx -\bb_t)$ be a self-concordant function. The following holds:
    \[\norm{\DD(\yy_1)\DD^{-1}(\yy_2)}_{\DD(\yy_1)}\le \frac{1}{(1-\norm{\yy_1-\yy_2}_{\DD(\yy_1)})^2}.\]
\end{lemma}

\begin{proof}
    The inequalities follow directly from the definition of self-concordance~\cite{renegar}. We have:
    \begin{align*}
        &\norm{\DD(\yy_1)\DD^{-1}(\yy_2)}_{\DD(\yy_1)} = \norm{\left(\DD(\yy_1)\DD^{-1}(\yy_2)\right)^{-1}}_{\DD(\yy_1)}\\
        &\quad\quad\quad\quad\quad=\min_{\mathbf{v}\ne 0}\frac{\norm{\mathbf{v}}_{\DD(\yy_1)}^2}{\mathbf{v}^\top\DD(\yy_1)\DD^{-1}(\yy_1)\DD(\yy_2)\mathbf{v}}\\
       &\quad\quad\quad\quad\quad \le \sup_{\mathbf{v}\ne 0} \frac{\norm{\mathbf{v}}^2_{\DD(\yy_2)}}{\norm{\mathbf{v}}^2_{\DD(\yy_1)}} \le \left(\frac{1}{1-\norm{\yy_2-\yy_1}_{\DD(\yy_1)}}\right)^2,
    \end{align*}
    where the last inequality follows from \eqref{def:selfconcordance}.
\end{proof}

%\begin{lemma}
 %   If $\exists m \in \mathbb{R}$ such that $\norm{\nabla^2\phi(\xx)^-1}\le \frac{1}{m}$, for all $\xx\in \Df$, then,
  %  \[\norm{\DD^{-1}(\yy)} \le \frac{1}{m}\quad \forall \yy \quad\text{\normalfont s.t. } \xx\in \Df.\]
%\end{lemma}
%\begin{proof}
 %   By construction we have:
 %   \begin{align*}
 %   \norm{\DD(\yy)\mathbf{v}} \ge &\norm{\nabla^2 \phi(\xx)\mathbf{v}} \ge m\norm{\mathbf{v}}\quad \forall \yy,\xx,
 %   \mathbf{v},\\
%\end{align*}
%and the result follows from~\cite[Lemma 1]{OnlineNewton}. This condition is equivalent to requiring $\phi$ to be strongly convex which will always be the case if $\Df$ is bounded.
%\end{proof}

The dual-feasible nature of every point on the central path provides an estimation of the error for each iteration. Defining the next iterate as $\yy^+ =\yy + n_t(\yy,\nn)$, Lemma~\ref{lem:nred} formalizes this concept. 
\begin{lemma}
\label{lem:nred}
    For any pair $\{\yy,\nn\}$ such that $\normy{n_t(\yy,\nn)}\le 1$, we have:
    \begin{equation}\label{normreduction}
\norm{n_t(\yy^+,\nn)}_{\DD(\yy^+)} \le \frac{\norm{n_t(\yy,\nn)}^2_{\DD(\yy)}}{\left(1-\norm{n_t(\yy,\nn)}_{\DD(\yy)}\right)^2}.\end{equation}
    
\end{lemma}

\textcolor{Blue}{The proof of this Lemma is provided in Appendix~A.}

An important consequence of Lemma \ref{lem:nred} is if we have an initial pair $\{\yy,\nn\}$ such that $\norm{n_t(\yy,\nn)}_{\DD(\yy)} \le \frac{1}{4}$, the next iterate satisfies the following inequality:
\[\norm{n_t(\yy^+,\nn)}_{\DD(\yy^+)}\le \frac{1}{9}.\]

It is possible to determine a point's proximity to the optimum using the norm of the Newton's step computed from this point. This result is presented in Lemma~\ref{lem:closeness}.

\begin{lemma}
\label{lem:closeness}
    For any pair $\{\yy, \nn\}$ such that $\normy{n_t(\yy,\nn)}\le \frac{1}{4}$, we have:
    \[\norm{\yy^\nn-\yy}_{\DD(\yy)} \le \normy{n_t(\yy,\nn)} + \frac{3\normy{n_t(\yy,\nn)}^2}{(1-\normy{n_t(\yy,\nn)})^3},\]
where $\yy^\nn=\begin{bmatrix}\xx^\nn\\ \vv^\nn\end{bmatrix}$ denotes the optimal primal and dual variables such that:
\begin{align*}
    \xx^\nn &\in \argmin_{\xx} \nn \cc^\top \xx-\phi(\xx)\\
    &\quad\quad\text{\normalfont s.t.}\quad\nn\cc-\nabla\phi(\xx^\nn)+\AA^\top\vv^\nn = 0.
\end{align*}
\end{lemma}

\begin{proof}
   \textcolor{Blue}{The proof follows the analysis from \cite[Section 2.2.5]{renegar} where Theorem 2.2.4 is substituted by Lemma \ref{lem:nred}.}
In the specific case where $\norm{n_t(\yy^+,\nn)}_{\DD(\yy^+)}\le \frac{1}{9}$, this implies:
\textcolor{Blue}{$\norm{\yy^\nn-\yy^+}_{\DD(\yy^+)} \le \frac{1}{6}$.}
\end{proof}

\begin{lemma}\label{lem:close2} If $\normy{\yy^\nn-\yy} \le \frac{1}{6}$, where $\yy$ is a feasible point, then: \textcolor{Blue}{$\norm{\yy^\nn-\yy}_{\DD(\yy^\nn)} \le \frac{1}{5}.$}
\end{lemma}
\begin{proof}
From the definition of self-concordance we have:
\begin{align*}
\frac{\norm{\yy^\nn-\yy}_{\DD(\yy^\nn)}}{\normy{\yy-\yy^\nn}} &\le \frac{1}{1-\normy{\yy-\yy^\nn}}\\
    \norm{\yy^\nn-\yy}_{\DD(\yy^\nn)} &\le \frac{\normy{\yy-\yy^\nn}}{1-\normy{\yy-\yy^\nn}}= \frac{\frac{1}{6}}{\frac{5}{6}}= \frac{1}{5}.
\end{align*}
\end{proof}

\begin{lemma}
\label{lem:b}
    Consider problem \eqref{prob2}. Let $\Delta\bb_t \stackrel{\text{def}}{=}\bb_{t+1}-\bb_t$ be the change in the online parameter $\bb_t$ between rounds.
    If $\norm{n_t(\yy,\nn)}\le\frac{1}{9}$ and $\norm{\Delta\bb_t}\le \sqrt{\frac{3m}{160}}$ for all $t$, then:\begin{align*}\norm{n_{t+1}(\yy,\nn)} \le \frac{1}{4}.\end{align*}
    
\end{lemma}

\begin{proof}
First, we remark that the only change between rounds is the parameter $\bb_t$. Therefore,
\[r_{t+1}(\yy,\nn)=r_t(\yy,\nn)+\begin{bmatrix}\mathbf{0}\\ \Delta\bb_t\end{bmatrix}.\]
Let $\Delta \mathbf{r}_t = \begin{bmatrix}\mathbf{0}\\ \Delta\bb\end{bmatrix}$. Thus, $\norm{\Delta \mathbf{r}_t}=\norm{\Delta \bb_t}$. The Hessian norm of the Newton's step becomes:
\begin{align*}
    \normy{n_{t+1}(\yy,\nn)}^2 &= (r_t(\yy,\nn)+\Delta \mathbf{r})^\top\DD(\yy)^{-1}(r_t(\yy,\nn)+\Delta \mathbf{r}_t)\\
    \le 2 &r_t(\yy,\nn)^\top\DD(\yy)^{-1}r_t(\yy,\nn) %\\&\quad\quad\quad\quad\quad\quad
    + 2 \Delta \mathbf{r}_t^\top \DD(\yy)^{-1} \Delta \mathbf{r}_t\\
    =& \frac{2}{81} + 2\norm{\Delta \bb_t}^2\norm{\DD(\yy)^{-1}}\\
    =& \frac{2}{81} + \frac{2\norm{\Delta \bb_t}^2}{m}\le \frac{1}{16},
\end{align*}
\textcolor{Blue}{where we used Assumption 4 and the bound on $\norm{\Delta\bb}$ to obtain the final inequality}.
\end{proof}
%
% The following is an interesting property of the Hessian norm of the constrained Newton step for feasible points.

% \begin{align*}
%     \norm{n(\yy)}_\yy^2 &= \begin{bmatrix}
%         \nabla d(\xx) +\AA^\top\vv\\ \mathbf{0}
%     \end{bmatrix}^\top\begin{bmatrix}\nabla^2 d(\xx) & \AA^\top\\ \AA& \mathbf{0}\end{bmatrix}^{-1}\begin{bmatrix} \nabla d(\xx) +\AA^\top\vv\\ \mathbf{0}\end{bmatrix}\\
%     &= \begin{bmatrix}
%         \nabla d(\xx) +\AA^\top\vv\\ \mathbf{0}
%     \end{bmatrix}^\top\begin{bmatrix}\ddi - \ddi\AAt(\AA\ddi\AAt)^{-1}\AA\ddi & \ddi\AAt(\AA\ddi\AAt)^{-1}\\
%     (\AA\ddi\AAt)^{-1}\AA\ddi & -(\AA\ddi\AAt)^{-1}\end{bmatrix}\begin{bmatrix}
%         \nabla d(\xx) +\AA^\top\vv\\ \mathbf{0}
%     \end{bmatrix}
% \end{align*}
% The structure of the problem is such that the following is true. Here, we use the fact that the inverse of any symmetric matrix is symmetric and our Hessian matrix is always symmetric.
% \begin{align*}
%     \norm{n(\yy)}_\yy^2 &= (\nabla d(\xx)^\top + \vv^\top\AA)(\ddi-\ddi\AAt(\AA\ddi\AAt)^{-1}\AA\ddi)(\nabla d(\xx)+\AA^\top\vv)\\
%     &= \nabla d(\xx)^\top\ddi\nabla d(\xx) - \nabla d(\xx)^\top \ddi\AAt(\AA\ddi\AAt)^{-1}\AA\ddi\nabla d(\xx)
% \end{align*}
% Although this formulation is quite heavy, the important element is that the Hessian norm of the Newton step is completely independent of the dual variable $\vv$ when the current iterate is feasible. Furthermore, we can prove that this is the same quantity as would be obtained from applying Newton's method to the equivalent projected problem. 

From~\cite[Section 2.3.1]{renegar}, we have that the restriction of any barrier functional to a subspace or translation of a subspace results in a barrier functional whose complexity barrier is smaller than that of the original functional. This implies that if $\norm{\phi(\xx)}_{\nabla^2 \phi(\xx)}^2\le v_f$ any pair $\xx, \vv$ satisfies:
\[\begin{bmatrix}
        \nabla \phi(\xx) +\AA^\top\vv\\ \mathbf{0}
    \end{bmatrix}^\top\DD(\xx)^{-1}\begin{bmatrix} \nabla \phi(\xx) +\AA^\top\vv\\ \mathbf{0}\end{bmatrix} \le v_f.\]
Defining  $h(\yy) \stackrel{\text{def}}{=} \DD(\xx)^{-1}\begin{bmatrix} \nabla \phi(\xx) +\AA^\top\vv\\ \mathbf{0}\end{bmatrix},$ we notice that $\normy{h(\yy)} \le \sqrt{v_f}$. We will use this result in Lemma \ref{lem:eta}.

\begin{lemma}
    \label{lem:eta}
    \textcolor{Blue}{For a pair $\{\yy,\nn\}$ such that $\yy$ is feasible and $\normy{n_t(\yy,\nn)}\le\frac{1}{9}$}, there exists a parameter $1<\beta\le 1+\frac{1}{8\sqrt{v_f}}$ such that for $\nn^+ = \beta\nn$, we have:
    \[\normy{n_t(\yy,\nn^+)} \le \frac{1}{4}.\]
\end{lemma}

\textcolor{Blue}{The proof is presented in Appendix~B. We now present a lemma linking optimal values to feasible points using the Hessian norm.}

\begin{lemma}
\label{lem:renmyst}
    Let $\xx^\nn_t\in \argmin_{\xx} \{ \nn \cc^\top\xx + \phi(\xx)$ {\normalfont s.t.} $\AA\xx-\bb_t = \mathbf{0}\}$. If $\bar{\xx}$ is a feasible point then,
    \begin{align}
\cc^\top (\bar{\xx}-\xx^\nn_t) &\le   \frac{v_f}{\nn_t}(1+\norm{\bar{\xx}-\xx^\nn_t}_{\nabla^2\phi(\xx^\nn_t)}).\end{align}
\end{lemma}
\begin{proof}
    The above identity follows from~\cite[Section 2.4.1]{renegar} for a linear objective.
\end{proof}

Finally, we relate $\norm{\yy_t-\yy_t^\nn}_{\DD(\yy_t)}$ to $\norm{\xx_t-\xx_t^\nn}_{\nabla^2\phi(\xx_t)}$ in Lemma~\ref{lem:yx}.
\begin{lemma}
    \label{lem:yx}
    For any primal-dual feasible point $\yy$, we have:
    \[\normy{\yy_t-\yy_t^\nn}\ge\norm{\xx_t-\xx_t^\nn}_{\nabla^2\phi(\xx_t)}.\]
\end{lemma}
\begin{proof}
    Recall from the definition of the Hessian norm:
    \begin{align*}
        \norm{\yy_t-\yy_t^\nn}^2_{\DD(\yy_t)} &= (\yy_t-\yy_t^\nn)^\top \DD(\yy_t)(\yy_t-\yy_t^\nn)\\
        &\ge \begin{bmatrix}
            \xx_t-\xx_t^\nn\\ \mathbf{0}
        \end{bmatrix}^\top \DD(\yy_t)\begin{bmatrix}\xx_t-\xx_t^\nn\\ \mathbf{0}\end{bmatrix}\\
        &= (\xx_t-\xx^\nn_t)^\top\nabla^2\phi(\xx_t)(\xx_t-\xx^\nn_t)\\
        &= \norm{\xx_t-\xx_t^\nn}_{\nabla^2\phi(\xx_t)}^2.
    \end{align*}
    Taking the square root on both sides, we obtain the desired result completing the proof.
\end{proof}

\section{Algorithms} \label{algos}

Suppose we have an initial primal-dual point $\xx_0$, $\vv_0$ and an initial barrier parameter $\nn_0$ such that $\norm{n_0(\yy_0, \nn_0)}_{\DD(\yy_0)} \le \frac{1}{9}$ and $\AA\xx_0=\bb_0$. \textcolor{Blue}{In practice, this can be achieved using offline optimization before the start of the online procedure.} We propose the online interior-point method for time-varying equality constraints (\texttt{OIPM-TEC}) presented in Algorithm~\ref{alg:TIP2}.

\begin{algorithm}[h!]
\caption{\texttt{OIPM-TEC}}\label{alg:TIP2}
\begin{algorithmic}
\State\textbf{Parameters:} $\AA$, $\cc$, $\beta$
\State\textbf{Initialization}: Receive $\xx_0 \in\mathbb{R}^n$, $\vv_0\in\mathbb{R}^p$, $\nn_0 \in\mathbb{R}^+$ such that $\norm{n_0(\yy_0, \nn_0)}_{\DD(\yy_0)}\le \frac{1}{9}$.
\FOR{$t=0,1,2...T$}
\State Implement the decision $\xx_t$.
\State Observe the outcome $\cc^\top\xx_t$ and the new constraint $\bb_t$.
\State \underline{Perform $t$-step:}
\State \hspace{-0.4cm}  $\begin{bmatrix}\tilde{\xx}_t\\ \tilde{\vv}_t\end{bmatrix}$ $=\begin{bmatrix}\xx_t\\ \vv_t\end{bmatrix}-\begin{bmatrix}\nabla^2 \phi(\xx_t) & \AA^\top\\ \AA& \mathbf{0}\end{bmatrix}^{-1}\begin{bmatrix} \nabla d_{\nn_t}(\xx_t) +\AA^\top\vv_t\\ \AA\xx_t-\bb_t\end{bmatrix}$.
\State \underline{Perform $\nn$-step:}
\State  $\nn_{t+1} = \nn_t\beta$,
\State \hspace{-0.5cm} $\begin{bmatrix}\xx_{t+1}\\ \vv_{t+1}\end{bmatrix}$ $=\begin{bmatrix}\tilde{\xx}_t\\ \tilde{\vv}_t\end{bmatrix}-\begin{bmatrix}\nabla^2 \phi(\tilde{\xx}_t) & \AA^\top\\ \AA& \mathbf{0}\end{bmatrix}^{-1}\begin{bmatrix} \nabla d_{\nn_{t+1}}(\tilde{\xx}_t) +\AA^\top\vv_t\\ \mathbf{0}\end{bmatrix}$.
\ENDFOR
\end{algorithmic}
\end{algorithm}

\begin{theorem}
    If Assumptions 1-4 hold, \textcolor{Blue}{$0 < \beta \le1+\frac{1}{8\sqrt{v_f}}$,  $\max_t \norm{\bb_t-\bb_{t-1}}\le \sqrt{\frac{3m}{160}}$}, and supposing an initial point $\xx_0$, $\vv_0$ and an initial barrier parameter $\nn_0$ such that $\norm{n_0(\yy_0, \nn_0)}_{\DD(\yy_0)} \le \frac{1}{9}$, then the regret and constraint violation of Algorithm \ref{alg:TIP2} are bounded by:
    \begin{align}
        R_{\text{{\normalfont d}}}(T)&\le \frac{11 v_f\beta}{5\nn_0(\beta-1)} + cV_T\\
        \text{{\normalfont Vio}}(T)&\le V_\bb.
    \end{align}
    The dynamic regret and constraint violation are, therefore, $O(V_T+1)$ and $O(V_\bb)$, respectively.
\end{theorem}

\begin{proof} We first derive some properties of the algorithm. 
Suppose at \textcolor{Blue}{some} round we have $\yy_{t-1}$ such that $\norm{n_{t-1}(\yy_{t-1}, \eta_{t-1})}_{\DD(\yy_{t-1})}\le \frac{1}{9}$. Then, when observing the new vector $\bb_t$, Lemma \ref{lem:b} yields $\norm{n_t(\yy_{t-1}, \eta_{t-1})}_{\DD(\yy_{t-1})}\le \frac{1}{4}$. By Lemma \ref{lem:nred}, this implies that the point $\tilde{\yy}$ obtained after the $t$-step will respect:
\[\norm{n_t(\tilde{\yy}_t, \eta_{t-1})}_{\DD(\tilde{\yy}_t)}\le \frac{1}{9}.\]
After updating $\nn$, Lemma \ref{lem:eta} yields:
\[\norm{n_t(\tilde{\yy}_t, \eta_t)}_{\DD(\tilde{\yy}_t)}\le \frac{1}{4}.\]
After the $\nn$-step from Lemma \ref{lem:nred} we have:
\[\norm{n_t(\yy_t, \eta_t)}_{\DD(\yy_t)}\le \frac{1}{9}.\]
Lemma \ref{lem:closeness} leads to:
\[\norm{\yy_t-\yy_{t-1}^{\nn_t}}_{\DD(\yy_t)}\le \frac{1}{6}.\]
Now, applying Lemma \ref{lem:close2}, we get:
\[\norm{\yy_t-\yy_{t-1}^{\nn_t}}_{\DD(\yy_{t-1}^{\nn_t})}\le \frac{1}{5}.\]
Finally, Lemma \ref{lem:yx} yields,
\begin{equation}\label{result}\norm{\xx_t-\xx_{t-1}^{\nn_t}}_{\nabla^2\phi(\xx_{t-1}^{\nn_t})}\le \frac{1}{5}.\end{equation}
It follows that because we initially have $\norm{n_0(\yy_0, \nn_0)}_{\DD(\yy_0)} \le \frac{1}{9}$ by assumption, then $\norm{n_{t-1}(\yy_{t-1}, \eta_{t-1})}_{\DD(\yy_{t-1})}\le \frac{1}{9}$ for all~$t$.

Defining $c=\norm{\cc}$, we now re-express the regret as:
\noindent \begin{align}\label{eq:regret2}
    R_{\text{d}}(T) &= \sum_{t=1}^{T} f(\xx_t)-f(\xx^*_t)\nonumber\\
    =& \sum_{t=1}^{T}\cc^\top\left(\xx_t - \xx^{\nn_t}_{t-1} +\xx^{\nn_t}_{t-1}-\xx^*_{t-1}+\xx^*_{t-1}-\xx^*_t \right)\nonumber\\
    \le& \sum_{t=1}^{T}\cc^\top\left(\xx_t - \xx^{\nn_t}_{t-1} +\xx^{\nn_t}_{t-1}-\xx^*_{t-1}\right)+c\norm{\xx^*_{t-1}-\xx^*_t }\nonumber\\
    =& \sum_{t=1}^T\cc^\top\big(\xx_t - \xx^{\nn_t}_{t-1}\big)+\sum_{t=1}^T\cc^\top\big(\xx^{\nn_t}_{t-1}-\xx^*_{t-1}\big)+cV_T.
\end{align}
Next, we bound the first terms of \eqref{eq:regret2}. The first sum represents the estimation error between the current decision and the optimum for a given $\nn$. Using Lemma \ref{lem:renmyst},
\begin{align*}
    \sum_{t=1}^T\cc^\top\big(\xx_t - \xx^{\nn_t}_{t-1}\big) &\le \sum_{t=1}^T \frac{v_f}{\nn_t}\left(1+\norm{\xx_t-\xx^{\nn_t}_{t-1}}_{\nabla^2 \phi(\xx^{\nn_t}_{t-1})} \right)\\
    &\le v_f\sum_{t=1}^T \frac{1}{\nn_t}\left(1+\frac{1}{5}\right)\\
    &= \frac{6v_f}{5\nn_0}\sum_{t=1}^T\left(\frac{1}{\beta}\right)^{t}\le \frac{6 v_f\beta}{5\nn_0(\beta-1)},
\end{align*}
where we used \eqref{result} to obtain the second inequality.

From~\cite[Chapter 9]{boyd}, we have that any optimizer $\xx_t^\nn$ of the functional $d_t(\xx_t, \nn)$ for a particular $\nn$ has an optimality gap of at most $\frac{v_f}{\nn}$, i.e., $\cc^\top(\xx_t^\nn-\xx_t^*)\le\frac{v_f}{\nn}$. Therefore,
\begin{align*}
    \sum_{t=1}^T\cc^\top(\xx^{\nn_t}_{t-1}-\xx^*_{t-1}) &\le  \sum_{t=1}^T \frac{v_f}{\nn_t}\\
    & = \frac{v_f}{\nn_0}\sum_{t=1}^T\left(\frac{1}{\beta} \right)^t\le \frac{v_f\beta}{\nn_0(\beta-1)}.
\end{align*}
Thus, the regret bound is : 
\[R_{\text{d}}(T) \le cV_T+\frac{11v_f\beta}{5\nn_0(\beta-1)}.\]

\textcolor{Blue}{
As for constraint violation, we have:
\begin{align*}
    \text{Vio}(T) = \sum_{t=1}^T \left[\sum_id_{K_i}(\xx_t)+\norm{\AA\xx_t-\bb_t}\right].
\end{align*}
We note that iterates are always in the feasible space defined by the inequality constraints $g_i(\xx)$. This is because, at each round, a unit-step is taken in the Newton direction which lies in the feasible space of the self-concordant functional~\cite[Section 2.2.1]{renegar}.} It follows that, $g_i(\xx_t)\preceq 0,\ \forall t,i$ and $d_{K_i}(\xx_t)=0$. We also observe that iterates are always feasible with respect to the previous round's equality constraints. This is again because of the unit-step taken during the $t$-step~\cite{boyd}. Using this fact, we can write:
\begin{align*}
    \text{Vio}(T) &= \sum_{t=1}^T\big[0 + \norm{\AA\xx_t-\bb_t}\big]\\
    &= \sum_{t=1}^T \norm{\bb_{t-1}-\bb_t}= V_{\bb},
\end{align*}
which completes the proof.
\end{proof}

When operating power grids, the exact optimal solution is not always required. Indeed there is a non-zero tolerance-region around the optimum that the electric grid operator can accept. For example, if the objective is cost-minimization, a solution to the OPF within a few cost units would fall within the tolerance of the grid operator. To translate this to a more traditional OCO framework, we can define an epsilon-regret $R_\epsilon(T)$ as:
\begin{equation}
    \label{def:epregret}
    R_\epsilon(T)=\sum_{t=1}^T \left[f_t(\xx_t)-f_t(\xx_t^*)-\epsilon\right]^+,
\end{equation} 
where $0<\epsilon $ is the desired tolerance. We remark that the epsilon-regret is null if all decisions are within our tolerance and positive otherwise. Note that epsilon-regret is a weaker metric than dynamic regret. A bounded dynamic regret implies a bounded epsilon-regret. \textcolor{Blue}{This metric removes the dependence of a neighbourhood term in the regret bound such as in~\cite{dallanese2020}.} 
\textcolor{Blue}{If instead of requiring sublinear dynamic regret, we are only interested in sublinear epsilon-regret, an algorithm requiring only a $t$-step can be used: the epsilon online interior-point method for time-varying equality constraints (\texttt{$\epsilon$OIPM-TEC}).} This algorithm is presented in Algorithm \ref{alg:EpsilonTIP}. 
\begin{algorithm}[h!]
\caption{\texttt{$\epsilon$OIPM-TEC}}\label{alg:EpsilonTIP}
\begin{algorithmic}
\State\textbf{Parameters:} $\AA$, $\cc$
\State\textbf{Initialization}: Receive $\xx_0 \in\mathbb{R}^n$, $\vv_0\in\mathbb{R}^p$, $\nn \in\mathbb{R}^+$ such that $\norm{n_0(\yy_0, \nn)}_{\DD(\yy_0)}\le \frac{1}{9}$.
\FOR{$t=0,1,2...T$}
\State Implement the decision $\xx_t$.
\State Observe the outcome $\cc^\top\xx_t$ and the new constraint $\bb_t$.
\State  \hspace{-0.5cm} $\begin{bmatrix}\xx_{t+1}\\ \vv_{t+1}\end{bmatrix}$ $=\begin{bmatrix}\xx_t\\ \vv_t\end{bmatrix}-\begin{bmatrix}\nabla^2 \phi(\xx_t) & \AA^\top\\ \AA& \mathbf{0}\end{bmatrix}^{-1}\begin{bmatrix} \nabla d_{\nn_t}(\xx_t) +\AA^\top\vv_t\\ \AA\xx-\bb_t\end{bmatrix}$.
\ENDFOR
\end{algorithmic}
\end{algorithm}

\begin{theorem}
    Let $\epsilon\ge 0$ be the tolerance. If Assumptions 1-4 hold, $0 < \beta \le1+\frac{1}{8\sqrt{v_f}}$, \textcolor{Blue}{$\max_t \norm{\bb_t-\bb_{t-1}}\le \sqrt{\frac{3m}{160}}$}, and we have the initial point $\xx_0$, $\vv_0$, and $\nn$ such that $\nn \ge \frac{11v_f}{5\epsilon}$, $\norm{n_0(\yy_0,\nn)}_{\DD(\yy_0)}<\frac{1}{9}$, then the epsilon-regret and constraint violation of Algorithm \ref{alg:EpsilonTIP} are bounded by:
    \begin{align}
        R_\epsilon(T) &\le cV_T\\
        \text{{\normalfont Vio}}(T) &\le V_\bb.
    \end{align}
\end{theorem}

\begin{proof}
    The constraint violation is bounded by $V_\bb$ following the same reasoning as for Algorithm \ref{alg:TIP2}.
    
    Similarly to Algorithm \ref{alg:TIP2}, the point $\yy_{t}$ obtained after the Newton's step will respect:
    $\norm{n_t(\yy_{t}, \eta)}_{\DD(\yy_t)}\le \frac{1}{9}.$
    Applying Lemmas \ref{lem:closeness}, \ref{lem:close2}, and \ref{lem:yx} successively, we have:
    $\norm{\xx_{t+1}-\xx_t^\nn}_{\nabla^2 \phi(\xx_t^\nn)}\le \frac{1}{5}$.
    We can therefore re-express the epsilon-regret as:
    \begin{align}
        R_\epsilon(T) &= \sum_{t=1}^T\left[f_t(\xx_t)-f_t(\xx_t^*)-\epsilon\right]^+\nonumber\\ 
    \end{align}
    \begin{align}
         %\le  \sum&_{t=1}^T\Big[\cc^\top\big(\xx_t - \xx^{\nn_t}_{t-1}+\xx^{\nn_t}_{t-1}-\xx^*_{t-1}%\nonumber\\&\quad\quad\quad\quad\quad\quad\quad\quad
        %+\xx^*_{t}-\xx^*_{t-1}\big)-\epsilon\Big]^+\nonumber\\
        \le \sum_{t=1}^T&\left[\frac{v_f}{\nn}\left(1+\norm{\xx_{t+1}-\xx_t^\nn}_{\nabla^2\phi(\xx_t^\nn)}\right)+\frac{v_f}{\nn}-\epsilon\right]^+ +cV_T\nonumber\\
        \label{eq:afs}
        \le \sum_{t=1}^T&\left[\frac{11v_f}{5\nn}-\epsilon\right]^++cV_T.      
    \end{align}
    Because $\nn \ge \frac{11v_f}{5\epsilon}$ by assumption, $\epsilon \ge \frac{11v_f}{5\nn}$. Therefore the leftmost sum of \eqref{eq:afs} is zero. This implies that $R_\epsilon(T)\le cV_T$, thus completing the proof.
\end{proof}

\section{Numerical Experiment}

We now illustrate the performance of our algorithms in an online OPF problem. In the online OPF, the grid operator seeks to serve time-varying loads within a fixed network while minimizing the generation costs. Consider a set of buses $\mathcal{B}\subset\mathbb{N}$ on which lies a set of generators $\mathcal{G}\subseteq\mathcal{B}$ connected to a set of loads $\mathcal{D}\subseteq \mathcal{B}$ by power lines $\mathcal{L}\subseteq \mathcal{B}\times\mathcal{B}$. The admittance of every line $ij\in\mathcal{L}$ is given by $\gamma_{ij}\in\mathbb{C}$.  

At time $t$, this problem's decision variables are the real and reactive power injections at every bus given by $p_{i,t}$ and $q_{i,t}$, respectively. The dependent variables are the voltages, $v_{i,t}$ at each bus. The parameters are the minimum and maximum power limits on active and reactive generation ($\underline{p},\underline{q},\overline{p},\overline{q})$, the apparent power limits in the power lines ($\overline{k}_{ij}$) and the voltage limits $(\underline{v}, \overline{v})$. The active and reactive demand at every bus at time $t$ is given by $p^d_{i,t}$ and $q^d_{i,t}$, respectively. Finally, in order to obtain a second-order cone relaxation of the problem, the matrix $\WW\in \mathbb{C}^{\text{card} \mathcal{B}^2}$ is formed where every element $\WW_{ij}$ for $ij\in\mathcal{L}$ represents the product $v_iv_j^*$, with $^*$ denoting the complex conjugate~\cite{taylor2015convex}. At time $t$, the online OPF problem takes the following form:
\begin{align}
    \min& \quad s_t\nonumber\\
        \text{s.t.}\quad &\sum_{i\in\mathcal{G}} a_i p_{i,t}^2+b_i p_{i,t}-s_t \le 0\nonumber\\ 
\label{OPF}
    & \underline{p_i} \le p_{i,t} \le \overline{p_i}\quad \forall i\in\mathcal{G}\nonumber\\
        & \underline{q_i} \le q_{i,t} \le \overline{q_i}\quad \forall i\in\mathcal{G}\\
        &\underline{v}^2 \le \WW_{ii,t}\le \overline{v}^2\quad\forall i\in\mathcal{B}\nonumber\\
        &\norm{(\WW_{ii,t}-\WW_{ij,t})\gamma^*_{ij}} \le \overline{k}_{ij}\quad \forall ij\in\mathcal{L}\nonumber\\
        &\norm{\begin{matrix}2\WW_{ij,t}\\ \WW_{ii,t}-\WW_{jj,t}\end{matrix}}\le \WW_{ii,t} + \WW_{jj,t}\quad \forall{ij\in\mathcal{L}}\nonumber\\
         p_{i,t} +&p^d_{i,t} + \jmath (q_{i,t}+q^d_{i,t}) = \sum_{j\in\mathcal{N}} (\WW_{ii,t}-\WW_{ij,t})\gamma^*_{ij}\quad \forall{ij\in\mathcal{L}},\nonumber
        %&p_{d,t} = q_{d,t} = 0 \quad\forall d\in\mathcal{G}\\
        %&p_{i,t} = q_{i,t} = 0 \quad \forall i \notin\mathcal{G},
\end{align}
where $\jmath$ is the imaginary number.
We consider a time-varying extension of the second-order cone relaxation of the Matpower 33-bus case network presented in~\cite{Matpower}. To represent demand fluctuation, a random variation $\Delta p_{d,t}$ is added to every load at every timestep. For this experiment, every variation is independent and taken as: $\Delta p_{d,t} = \frac{0.01\zeta}{\sqrt{t}}$ where $\zeta$ is uniformly sampled in $[0,1]$ at every round. The load variation diminishes with the square root of $t$ ensuring that the variation in optima ($V_T$) and in the online parameter $\bb_t$ ($V_\bb$) are sublinear. If the offline problem, at a given timestep, is unsolvable, new random variables $\zeta$ are selected. The time horizon is set at $T=2000$.
For \texttt{OIPM-TEC}, the initial barrier parameter is set to $\eta_0=1$ and the parameter $\beta$ is set to $1.02$. For \texttt{$\epsilon$OIPM-TEC}, the barrier parameter is set to $\eta = 10^4$. \textcolor{Blue}{This value was chosen as it respects the condition $\eta \ge \frac{11 v_f}{5\epsilon}$ for the selected~$\epsilon$.} The initial primal-dual point $\xx_0, \vv_0$ for both algorithms is determined from the offline optimal solution of the Matpower 33-bus network. The acceptable tolerance $\epsilon=0.015$ is selected for the epsilon-regret which represents less than $0.001\%$ of the offline optimal cost. This value is equivalent to $0.015\$$/hr in this case. As a comparison, the \texttt{MOSP} algorithm~\cite{MOSP} is also applied to problem \eqref{OPF} using the same initial primal point as \texttt{OIPM-TEC}. In the \texttt{MOSP} implementation, the time-varying equality constraints are relaxed to inequality constraints such that $\AA\xx_t-\bb_t\le\mathbf{0}$. This enables the operator to over-serve loads and should be strictly respected at the optimum. The other linear inequality constraints are implemented directly while a projection is used to satisfy the second-order cone constraints. The parameters of the algorithm are set to : $\alpha=\mu=t^{-\frac{1}{3}}$, where $t$ denotes the current round.

\begin{figure}
    \centering
        \includegraphics[width=0.45\textwidth]{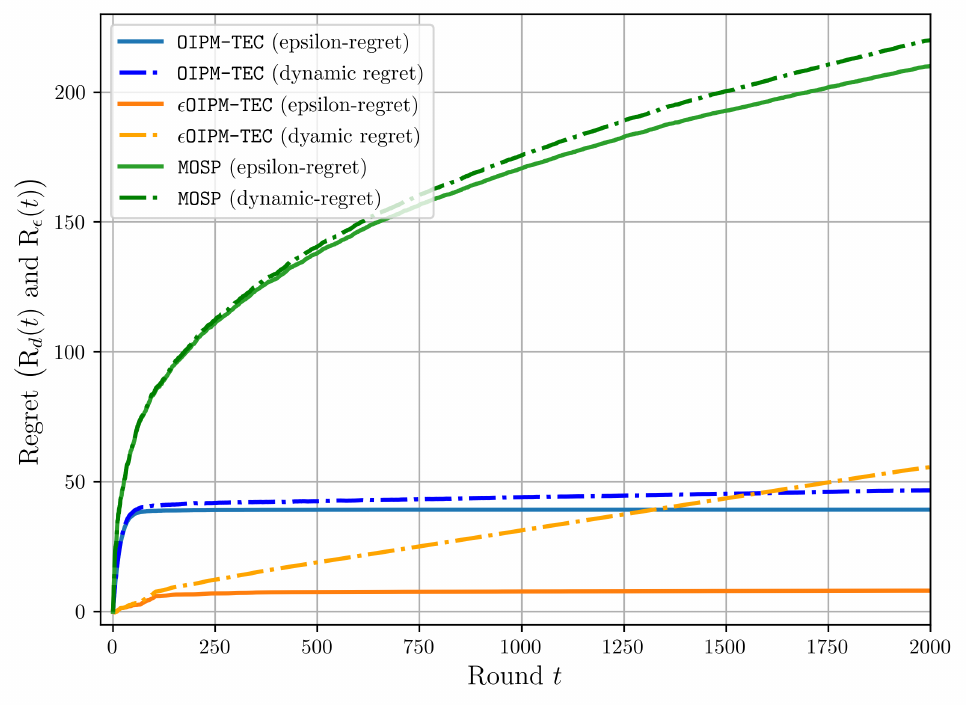}
        \vspace{-0.5cm}
        \caption{Regret comparison of \texttt{OIPM-TEC}, \texttt{$\epsilon$OIPM-TEC}, and \texttt{MOSP}}
        \label{fig:regret2}
\end{figure}

\begin{figure}
    \centering
    \includegraphics[width=0.45\textwidth]{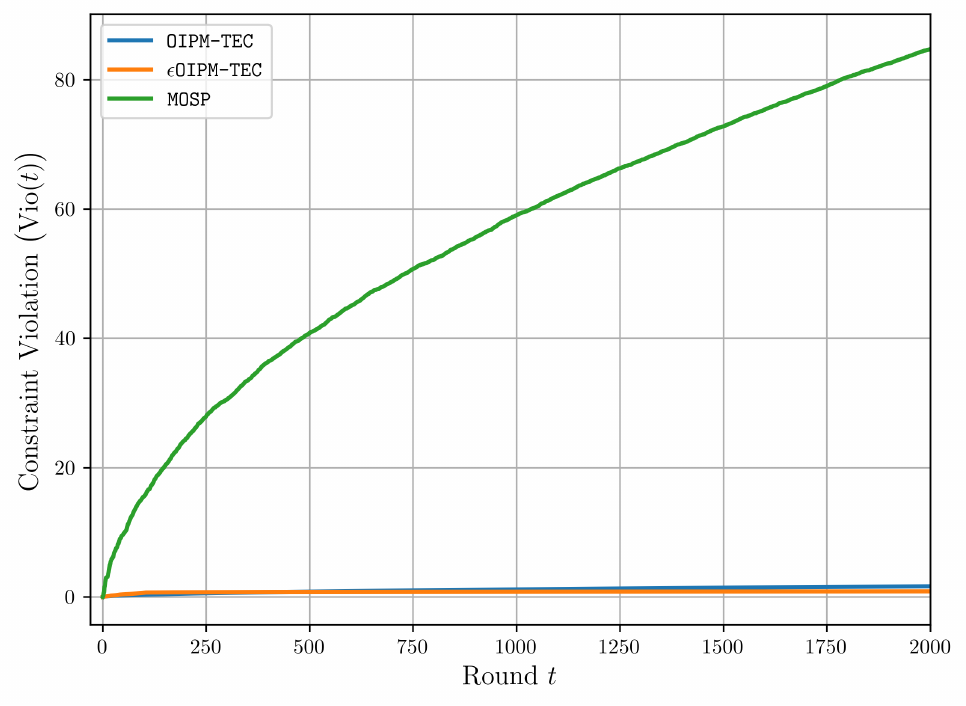}
    \vspace{-0.5cm}
    \caption{Constraint violation comparison of \texttt{OIPM-TEC}, \texttt{$\epsilon$OIPM-TEC}, and \texttt{MOSP}}
    \label{fig:const2}
    \vspace{-0.5cm}
\end{figure}

The resulting regret (dynamic regret and epsilon regret are considered) and constraint violation suffered by the three algorithms are presented in Figures \ref{fig:regret2} and \ref{fig:const2}. We remark that all algorithms achieve sublinear epsilon-regret and constraint violation. However, the epsilon-regret suffered by \texttt{OIPM-TEC}, while higher than that of \texttt{$\epsilon$OIPM-TEC}, is much lower than for \texttt{MOSP}. Both \texttt{OIPM-TEC} and \texttt{$\epsilon$OIPM-TEC} exhibit sublinear constraint violation terms that outperform \texttt{MOSP}. \textcolor{Blue}{Another element to consider is that $\epsilon$\texttt{OIPM-TEC} initially outperforms \texttt{OIPM-TEC} in terms of dynamic regret. This is due to the initally higher $\nn$ value of the former. }

\section{Conclusion}
In this work, an online optimization algorithm using an interior-point method-like update capable of solving the online OPF is presented. This algorithm, the \texttt{OIPM-TEC}, admits time-varying equality constraints and generalized inequalities while possessing a $O(V_T+1)$ dynamic regret bound and a $O(V_\bb)$ constraint violation bound. Additionally, if there exists an acceptable tolerance when solving an online optimization problem, the \texttt{$\epsilon$OIPM-TEC} algorithm is capable of ensuring the tightest possible bound with respect to $T$ on regret and constraint violation. These two algorithms are applied to an online OPF example showcasing their performance compared to other methods from the literature. 

The possibility of using OCO approaches to solve the OPF bridges the gap between traditional OCO algorithms with performance guarantees and real-time OPF approaches like~\cite{realtimeOPF} and~\cite{OptimalPowerFlowPursuit} that enjoy good real-world performance. A possible future direction is the development of online primal-dual interior-point methods which could increase numerical stability. Performance of online interior-point for fully time-varying problems is also of interest. Another avenue is considering inexact inversion techniques such as in~\cite{onlinebfgs}. 

\section*{References}
\bibliographystyle{IEEEtran}
\bibliography{ref.bib}

\newpage
\appendix
\subsection{Proof of Lemma \ref{lem:nred}}

From the definition of the Newton's step, we have:
\[\norm{n_t(\yy, \nn)}_{\DD(\yy)}^2 = r_t(\yy,\nn)^\top\DD(\yy)^{-1}r_t(\yy,\nn).\]
We can express $\norm{n_t(\yy^+,\nn)}_{\DD(\yy^+)}$ as:
\begin{align*}
    &\norm{n_t(\yy^+,\nn)}^2_{\DD(\yy^+)} = r_t(\yy^+,\nn)^\top\DD(\yy^+)^{-1}r_t(\yy^+,\nn)\\
    &\quad\quad= r_t(\yy^+,\nn)^\top\DD(\yy)^{-1}\DD(\yy)\DD(\yy^+)^{-1}r_t(\yy^+,\nn)\\
    &\quad\quad\le \norm{\DD(\yy)\DD(\yy^+)^{-1}}_{\DD(\yy)}r_t(\yy^+,\nn)^\top\DD(\yy)^{-1}r_t(\yy^+,\nn)\\
    &\quad\quad= \norm{\DD(\yy)\DD(\yy^+)^{-1}}_{\DD(\yy)} \norm{\DD(\yy)^{-1}r_t(\yy^+,\nn)}^2_{\DD(\yy)},
\end{align*}
where $\DD(\yy)^\top = \DD(\yy)$ because $\DD$ is symmetric.
Using Assumption 4, we can bound the first term by:
\begin{align*}
    \norm{\DD(\yy)\DD(\yy^+)^{-1}}_{\DD(\yy)} &\le \frac{1}{(1-\norm{\yy^+-\yy}_{\DD(\yy)})^2}\\
    &= \frac{1}{(1-\norm{n_t(\yy,\nn)}_{\DD(\yy)})^2.}
\end{align*}
Therefore, we have that:
\begin{equation}\label{eq:22}\norm{n_t(\yy^+,\nn)}_{\DD(\yy^+)}\le \frac{\norm{\DD(\yy)^{-1}r_t(\yy^+,\nn)}_{\DD(\yy)}}{1-\norm{n_t(\yy,\nn)}_{\DD(\yy)}}.\end{equation}
We now bound the numerator of \eqref{eq:22}. We have:
\begin{align*}
    &\normy{\DD(\yy)^{-1}r_t(\yy^+,\nn)}\\
    &\quad\quad=\normy{\DD(\yy)^{-1}r_t(\yy^+,\nn)-\DD(\yy)^{-1}r_t(\yy,\nn)+n_t(\yy,\nn)}\nonumber\\ %\label{eq1}
    &\quad\quad\le\normy{\DD(\yy)^{-1}\big[r_t(\yy^+, \eta)-r_t(\yy, \eta)\big]}+ \normy{n_t(\yy,\nn)}\\
    &\quad\le \int_0^1 \normy{\DD(\yy)^{-1}\DD\big(\yy-\tau n_t(\yy,\nn)\big)}\text{d}\tau \nonumber+ \int_0^1\normy{n_t(\yy,\nn)}\text{d}\tau,\end{align*} Where we use the fundamental theorem of calculus to obtain the last inequality. Using Lemma \ref{invHess} yields: \begin{align} %\label{eq2}
    &\normy{\DD(\yy)^{-1}r_t(\yy^+,\nn)}\nonumber\\
    &\quad\quad\le \int_0^1 \frac{\normy{n_t(\yy,\nn)}}{(1-\tau\normy{n_t(\yy,\nn)})^2}\text{d}\tau
    + \normy{n_t(\yy,\nn)}\nonumber\\
    &\quad\quad= \left[\frac{-1}{1-\tau\normy{n_t(\yy,\nn)}}\right]_0^1+\normy{n_t(\yy,\nn)}\nonumber\\
    &\quad\quad= \frac{\normy{n_t(\yy,\nn)}^2}{(1-\normy{n_t(\yy,\nn)})},\label{eq:3}
\end{align}
Substituting \eqref{eq:3} into \eqref{eq:22} completes the proof. $\blacksquare$

\subsection{Proof of Lemma \ref{lem:eta}}

We can decompose the Newton's step taken at any feasible point in the following way: 
\begin{align*}n_t(\yy, \nn) &= -\DD(\xx)^{-1}\begin{bmatrix}\nn \cc + \nabla\phi(\xx) + \AA^\top\vv\\ \mathbf{0}\end{bmatrix}\\
&= -\DD(\xx)^{-1}\left(\begin{bmatrix}\nn \cc \\ \mathbf{0}\end{bmatrix}+\begin{bmatrix}\nabla\phi(\xx) + \AA^\top\vv \\ \mathbf{0}\end{bmatrix}\right). \end{align*}
This allows us to express $n_t(\yy, \nn^+)$ as a function of $n_t(\yy, \nn)$:
\[n(\yy,\nn^+) = \frac{\nn^+}{\nn}n_t(\yy,\nn) + \left(\frac{\nn^+}{\nn}-1\right)h(\yy).\]
Taking the Hessian norm on both sides leads to:
\begin{align*}
    \normy{n_t(\yy,\nn^+)} &\le \frac{\nn^+}{\nn}\norm{n_t(\yy,\nn)}_{\DD(\yy)}+\left |1-\frac{\nn^+}{\nn}\right |\normy{h(\yy)}\\
    & \le  \frac{\nn^+}{\nn}\normy{n_t(\yy,\nn)}+\left |1-\frac{\nn^+}{\nn}\right |\sqrt{v_f}.
\end{align*}
By assumption, $\normy{n_t(\yy,\nn)}\le\frac{1}{9}$. Letting $\frac{\nn^+}{\nn}=\beta\le1+\frac{1}{8\sqrt{v_f}}$ we obtain:
\begin{align*}
    \normy{n_t(\yy,\nn^+)} &\le \frac{1}{9}\left(1+\frac{1}{8\sqrt{v_f}}\right)+\frac{\sqrt{v_f}}{8\sqrt{v_f}}\le \frac{1}{4},
\end{align*}
where we use the fact that $\sqrt{v_f}\ge 1$~\cite[Section 2.2]{renegar} to complete the proof. $\blacksquare$

\end{document}